\documentclass[12pt]{article}
\usepackage{latexsym}
\usepackage{amsmath,amsthm,amssymb}
\makeatletter
\newfont{\footsc}{cmcsc10 at 8truept}
\newfont{\footbf}{cmbx10 at 8truept}
\newfont{\footrm}{cmr10 at 10truept}
\renewcommand{\ps@plain}{%
\renewcommand{\@oddfoot}{\footsc the electronic journal of combinatorics
  {\footbf 10} (2003), \#R00\hfil\footrm\thepage}}
\makeatother
\newcommand{\bfrac}[2]{\left(\frac{#1}{#2}\right)}
\newcommand{\bin}[2]{\binom{#1}{#2}}

\newcommand{\brac}[1]{\left(#1\right)}
\def\cF{{\cal F}}
\def\cD{{\cal D}}

\def\cB{{\cal B}}

\def\cC{{\cal C}}

\def\cE{{\cal E}}
\def\cA{{\cal A}}

\def\om{\omega}

\def\Pr{\mbox{{\bf P}}}

\def\whp{{\bf whp}}

\newtheorem{lemma}{Lemma}
\newtheorem{theorem}{Theorem}

\newcommand{\proofend}{\hfill\mbox{$\Box$}}

\newcommand{\rdown}[1]{{\mbox{$ \lfloor #1 \rfloor $}}}

\newtheorem{remark}{Remark}
\newtheorem{claim}{Claim}

\def\Pr{\mbox{{\bf Pr}}}

\begin{document}
\pagestyle{plain}
\large
\date{}
\makeatletter \title{On Randomly Generated Intersecting Hypergraphs}
\author{Tom Bohman\thanks{ Supported in part by NSF grant DMS-0100400.}\\
\small Carnegie Mellon University,\\[-0.8ex]
\small \texttt{tbohman@andrew.cmu.edu}
\and Colin Cooper\\
\small Kings College,London\\[-0.8ex]
\small \texttt{ccooper@dcs.kcl.ac.uk}
\and Alan Frieze\thanks{Supported in part by NSF grant CCR-9818411.} \\
\small Carnegie Mellon University,\\[-0.8ex]
\small \texttt{alan@random.math.cmu.edu}
\and Ryan Martin\thanks{Supported in part by NSF VIGRE grant DMS-9819950.} \\
\small Iowa State University,\\[-0.8ex]
\small \texttt{rymartin@math.iastate.edu}
\and Mikl\'os Ruszink\'o\thanks{ Research was
partially supported by OTKA Grants T 030059, T 029074, T 038198, by the
Bolyai Foundation and by KFKI-ISYS Datamining Lab. and by a NATO
Collaborative Linkage Grant}\\[-0.8ex]
\small Computer and Automation Institute\\[-0.8ex]
\small Hungarian Academy of Sciences,\\[-0.8ex]
\small \texttt{ruszinko@lutra.sztaki.hu}}
\maketitle

\small Submitted: December 16, 2002;  Accepted: July 30, 2003.\\
\small MR Subject Classifications: 05D05, 05D40
\makeatother

\begin{abstract}
Let $c$ be a positive constant. We show that if $r=\rdown{cn^{1/3}}$ and the members of $\binom{[n]}{r}$ are
   chosen sequentially at random to form an intersecting hypergraph then
with limiting probability $(1+c^3)^{-1}$, as $n\to\infty$, the resulting family will be of maximum size
$\binom{n-1}{r-1}$.
\end{abstract}

\section{Introduction}

An {\em intersecting hypergraph} is one in which each pair of edges
has a non-empty intersection.  Here, we consider {\em $r$-uniform
hypergraphs} which are those for which all edges contain $r$
vertices.

The motivating idea for this paper is the classical
Erd\H{o}s-Ko-Rado theorem \cite{EKR} which states that a maximum size
$r$-uniform intersecting hypergraph  has $\binom{n-1}{r-1}$ edges if
$r\leq n/2$ and $\binom{n}{r}$ edges if $r>n/2$.  Furthermore, for
$r<n/2$ any maximum-sized family must have the property that all edges
contain a common vertex.

In the last four decades this theorem has attracted the attention of
many researchers and it has been generalized in many ways. It is worth
mentioning for example the famous conjecture of Frankl on the structure of
maximum $t$-intersecting families in a certain range of $n(t,r)$ which was
investigated by Frankl and F\"uredi \cite{FF} and completely solved only a
few years ago by Ahlswede and Khachatrian \cite{AK}. Another type of
generalization can be found in \cite{BFRT}.

The first attempt (and as far as we know the only one) to `randomize'
this topic was given by Fishburn, Frankl, Freed, Lagarias and Odlyzko
\cite{FFF}.  Also note that other random hypergraph structures were
considered already by R\'enyi e.g., in~\cite{R}, he identified the
anti-chain threshold.  Here we try to continue this line of
investigation. Our goal is to describe the structure of random
intersecting systems. More precisely, we consider taking edges
on-line; that is, one at a time, ensuring that at each stage, the
resulting hypergraph remains intersecting.  I.e., we consider the
following random process:

{\sc choose random intersecting system}\\
Choose $e_1\in \binom{[n]}{r}$.  Given
$\cF_i:=\{e_1,\ldots,e_i\}$, let $\cA(\cF_i)=\{e\in \binom{[n]}{r}:\;e\notin \cF_i$
and $e\cap e_j\neq\emptyset$ for $1\leq j\leq i\}$.
Choose $e_{i+1}$ uniformly at random from
$\cA(\cF_i)$. The procedure halts when $\cA(\cF_i)=\emptyset$ and $\cF=\cF_i$ is then
output by the procedure.

It should be made clear that sets are chosen {\em without} replacement.
\section{Definitions}

Let $[n]$ be the set of vertices of the hypergraph ${\cal H}$.

A {\em star} is collection of sets such that any pair in the collection has the same one-element intersection $\{x\}$,
which is referred to as the {\em kernel}.  A star with $i\geq 2$ edges is referred to as an $i$-star. A single
edge is a $1$-star, by convention.  We say that $\cal H$ is {\em fixed} by $x$ if every
member of $\cal H$ contains $x$.

For any sequence of events ${\cal E}_n$, we
will say that ${\cal E}_n$ occurs with high probability (i.e., \whp) if
$\lim_{n\to\infty}\Pr({\cal E}_n)=1$.

\section{The Erd\H{o}s-Ko-Rado Threshold}

The following theorem determines the threshold for the event
that edges chosen online to form an intersecting hypergraph will
attain the Erd\H{o}s-Ko-Rado bound.
\begin{theorem}\label{main}\
Let $\cE_{n,r}$ be the event that $|\cF|=\binom{n-1}{r-1}.$
For $r<n/2$, this is
equivalent to $\cal F$ fixing some $x\in [n]$.  Then if
$r=c_nn^{1/3}<n/2$,
$$\lim_{n\to\infty}\Pr(\cE_{n,r})=\begin{cases}1&c_n\to 0\\
\frac{1}{1+c^3}&c_n\to c\\0&c_n\to\infty\end{cases}.$$
\end{theorem}

\noindent{\bf Note:} If $r>n/2$, then all of $\binom{[n]}{r}$ is an
intersecting hypergraph.  If $r=n/2$ then for any ${\cal H}$ chosen
online to be an intersecting hypergraph, it will have size
$$\binom{n-1}{n/2-1}=\frac{1}{2}\binom{n}{n/2}.$$  In the case of $r=n/2$,
however, a vertex will not necessarily be fixed for even $n\geq 4$.

\section{Proof of Theorem~\ref{main}}
\subsection{Main Lemmas}
Before we prove relevant lemmas, we need to define some events.
\begin{itemize}
   \item Let $\cA_i$ be the event that $\cF_i$ is an $i$-star, for
   $i\geq 1$.
   \item Let $\cB_i$ be the event that $\cap_{j=1}^i
   e_j\neq\emptyset $, for $i\geq 3$.
   \item Let $\cC$ be the event that $e_3$ contains all of $e_1\cap
   e_2$ as well as at least one vertex in $(e_1\setminus e_2)\cup (e_2\setminus e_1)$.
   \item Let $\cD$ be the event that there is some $r$-set that
   intersects all currently chosen edges but fails to contain any vertex in their common
   intersection.
\end{itemize}

\begin{lemma}
If $r=o(n^{1/2})$ then
$$ \Pr(\cA_2)=1-o(1). $$
\label{propA2}
\end{lemma}

The fulcrum on which Theorem~\ref{main} rests is
Lemma~\ref{propA3}.
\begin{lemma}
If $r=o(n^{1/2})$ then
$$ \Pr(\cA_3)=\frac{1-o(1)}{1+\frac{(r-1)^3}{n}(1+o(1))} $$
\label{propA3}
\end{lemma}

\begin{lemma}
If $r=o(n^{2/5})$ and $m=O(n^{1/2}/r)$ then
$$ \Pr(\cA_{m}\mid\cA_3)=\exp\left\{-\frac{m^2r^2}{4n}+o(1)\right\} $$
\label{propAm+1}
\end{lemma}

\begin{remark}
Observe that Lemmas \ref{propA2}, \ref{propA3}, \ref{propAm+1} imply that if
$r=d_nn^{1/4}$, then the probability of the event $\cA_{r+1}$
approaches $\exp\{-d^4/4\}$ as $d_n\rightarrow d$. Furthermore, the
occurrence of $\cA_{r+1}$ immediately implies $\cA_s$ for $s>r+1$.
\end{remark}

\begin{lemma}
If $r=o(n^{1/2})$ then
$$ \Pr(\cC\mid \cA_2)=o(1). $$
\label{propC}
\end{lemma}

\begin{lemma}
If $r=o(n^{3/8})$ then
$$ \Pr(\cB_{3r}\mid\cA_4)=1-o(1) . $$
\label{propB3r}
\end{lemma}

\begin{lemma}
If $r=o(n^{2/5})$ then
$$ \Pr(\cD\mid\cB_{3r},\cA_4)=o(1). $$
\label{propD}
\end{lemma}

\begin{lemma}
If $r=\omega(n^{1/3})$ (i.e. $r/n^{1/3}\to\infty$) and $r=o(n^{2/3})$ then
$$ \Pr(\cB_3)=o(1). $$
\label{propB3}
\end{lemma}

\begin{lemma}
If $r=\om(n^{1/2})$ and $2\log_2n\leq m=o(e^{r^2/n})$ then
$$ \Pr(\cB_m)=o(1). $$
\label{propBm}
\end{lemma}
\subsection{Using these lemmas}
\noindent {\bf Case 1:} $r\leq n^{1/3}\log n$.

Suppose first that $c_n\to c$. Then
Lemma \ref{propA2} shows that $\cA_2$ occurs \whp.
Given $\cA_2$ there are 3 disjoint possibilities
\begin{equation}\label{three}
\cA_3\;\dot{\cup}\;\overline{\cB_3}\;\dot{\cup}\;\cC.
\end{equation}
 Lemma \ref{propC} shows that the conditional probability of $\cC$ tends to zero.
Lemma \ref{propA3} shows that
$\cA_3$ occurs with limiting probability $\frac{1}{1+c^3}$ and so given $\cA_2$
the probability of $\overline{\cB_3}$ tends to $\frac{c^3}{1+c^3}$.
If $\cB_3$ does not occur then $\cF$ cannot fix an element.

Suppose then that $\cA_3$ occurs and $e_1\cap e_2\cap e_3=\{v\}$.
We
use Lemma~\ref{propAm+1} with $m=4$ to show that $\cA_4$ occurs with conditional
probability $1-o(1)$. Then, given $\cA_4$ we can use
Lemma~\ref{propB3r} to show that $\cB_{3r}$ occurs \whp\ and Lemma \ref{propD}
to show that with conditional probability $1-o(1)$, $\cF$ must fix $v$.

If $c_n\to 0$ then $\cA_3$ occurs \whp\ and we conclude
as in the previous paragraph that with conditional probability $1-o(1)$, $\cF$ must fix $v$,
where $e_1\cap e_2\cap e_3=\{v\}$.

Now assume that $c_n\to\infty$.  We still have $\cA_2$ occuring \whp, but
now $\overline{\cA_3}$ occurs \whp. Using decomposition (\ref{three}) and Lemma \ref{propC}
to rule out event $\cC$ we see that $\overline{\cB_3}$ occurs \whp\ and so $\cF$ cannot fix any element.

\noindent {\bf Case 2:} $n^{1/3}\log n\leq r\leq  n^{1/2}\log n$. \\
Here we use Lemma~\ref{propB3}, which immediately gives that \whp\
$\cF_3$ has no vertex of degree 3; thus $\cF$ cannot fix any element.

\noindent {\bf Case 3:} $n^{1/2}\log n\leq r<n/2$. \\ In this case, we
apply Lemma~\ref{propBm} with
$m=\exp\left\{\frac{r^2}{3n}\right\}$ and we see that
$$ \Pr(\cB_m)= O\left(\exp\left\{-\frac{r^2}{3n}\right\}\right)=o(1) . $$

So $\cF_m$ fails, \whp, to have a vertex of degree $m$, in which case
$\cF$ cannot fix any element.  \proofend

\section{Proofs of Lemmas}

\subsection{Proof of Lemma~\ref{propA2}}

First we see that
\begin{eqnarray}
\Pr(\cA_1)&=&1.\label{1}\\
\Pr\left(\cA_{2}\mid \cA_1\right)&=&
\frac{r\binom{n-r}{r-1}}{\binom{n}{r}-\binom{n-r}{r}}\label{2}\\
&=&\frac{\frac{rn^{r-1}}{(r-1)!}\brac{1+O\bfrac{r^2}{n}}}
{\frac{n^r}{r!}\brac{1-1+\frac{r^2}{n}+O\bfrac{r^3}{n^2}}}\nonumber\\
&=&1+O\bfrac{r^2}{n}.\nonumber
\end{eqnarray}
\proofend

\subsection{Proof of Lemma~\ref{propA3}}

Continuing as in (\ref{2}),
\begin{equation}
\Pr\left(\cA_{i+1}\mid \cA_i\right)=
\frac{\binom{n-i(r-1)-1}{r-1}}{\binom{n-1}{r-1}+N_i-i},\qquad i\geq 2.
\label{3}
\end{equation}
For $i\geq 2$, the quantity $N_i$ is the number of $r$ sets that
intersect all of ${\cal F}_i$ but fail to contain the one-vertex
kernel of ${\cal F}_i$.  Thus,
\begin{equation}\label{44}
(r-1)^i\binom{n-i(r-1)-1}{r-i} \leq N_i\leq
(r-1)^i\binom{n-i-1}{r-i}.
\end{equation}
The lower bound comes from taking a single vertex (not the kernel) from
each of the edges and $r-i$ vertices from the remainder of the vertex set.
The upper bound comes from taking one vertex (not the kernel) from each of
the edges and $r-i$ other non-kernel vertices.

Simple computations give, for $r=o(n^{1/2})$,
\begin{eqnarray}
N_2&=&(1+o(1))\frac{(r-1)^3}{n}\binom{n-1}{r-1}. \label{+1}\\
N_3&\leq&(1+o(1))\binom{n-1}{r-1}.\label{4}\\
\binom{n-i(r-1)-1}{r-1}&=&(1+o(1))\binom{n-1}{r-1}.\label{4a}
\end{eqnarray}

It follows from (\ref{3}), (\ref{+1}), (\ref{4}) and (\ref{4a}) that
$$\Pr(\cA_3\mid \cA_2)=\frac{1-o(1)}{1+\frac{(r-1)^3}{n}(1+o(1))}.$$
Lemma~\ref{propA2} then gives that
\begin{equation}\label{A3}
\Pr(\cA_3)=\frac{1-o(1)}{1+\frac{(r-1)^3}{n}(1+o(1))}.
\end{equation}
\proofend

\subsection{Proof of Lemma~\ref{propAm+1}}

We estimate for $3\leq i\leq r$:
\begin{equation}\label{leq2}
\frac{(r-1)^i\binom{n-i-1}{r-i}}{\binom{n-1}{r-1}}\leq
\frac{r^i\binom{n-1}{r-i}}{\binom{n-1}{r-1}}
=O\brac{\frac{r^{2i-1}}{n^{i-1}}}.
\end{equation}
It then follows from (\ref{3}), (\ref{44}) and (\ref{leq2})
that for $3\leq i\leq r$,
\begin{eqnarray}
\Pr(\cA_{i+1}\mid \cA_i)&=&\frac{\binom{n-i(r-1)-1}{r-1}}
{\binom{n-1}{r-1}\brac{1+O\brac{\frac{r^{2i-1}}{n^{i-1}}}}}\nonumber\\
&=&1-\frac{ir^2}{2n}+O\brac{\frac{i^2r^3}{n^2}+\frac{r^{2i-1}}{n^{i-1}}}. \label{=2}
\end{eqnarray}
Equation (\ref{=2}) implies that
\begin{eqnarray*}
\Pr(\cA_{m+1}\mid \cA_3)&=&\prod_{i=3}^m\Pr(\cA_{i+1}\mid \cA_i)\\
&=&\prod_{i=3}^m\brac{1-\frac{ir^2}{2n}+O\brac{\frac{i^2r^3}{n^2}+\frac{r^{2i-1}}{n^{i-1}}}}\\
&=&\prod_{i=3}^m\exp\left\{-\frac{ir^2}{2n}+O\brac{\frac{i^2r^4}{n^2}+\frac{r^{2i-1}}{n^{i-1}}}\right\}\\
&=&\exp\left\{-\frac{m^2r^2}{4n}+o(1)\right\}.
\end{eqnarray*}
\proofend
\subsection{Proof of Lemma~\ref{propC}}

A simple computation suffices:
$$\Pr(\cC\mid\cA_2)\leq\frac{2r\bin{n-2}{r-2}}{\bin{n-1}{r-1}-2}
\leq\frac{2r^2}{n-2r\bin{n-1}{r-1}^{-1}}=O\left(\frac{r^2}{n}\right).$$
\proofend

\subsection{Proof of Lemma~\ref{propB3r}}

Assuming that both $\cA_4$ and $\cB_i$ occur for $i\geq 4$, there are
at most $(r-1)^4\binom{n-1}{r-4}$ $r$-sets which do not contain $v$
and which meet $e_1,e_2,e_3,e_4$.  On the other hand there are
$\binom{n-1}{r-1}-i$ $r$-sets which contain $v$ and are not edges of
$\cF_i$.  As a result, for $i\geq 4$,
\begin{equation}\label{10}
\Pr(\overline{\cB_{i+1}}\mid \cB_i,\cA_4)\leq
\frac{(r-1)^4\binom{n-1}{r-4}}{\binom{n-1}{r-1}-i}\leq
\frac{2r^7}{n^3}.
\end{equation}

Thus
\begin{eqnarray*}
\Pr(\cB_{3r}\mid \cA_4)&=&\prod_{i=4}^{3r-1}
\Pr(\cB_{i+1}\mid \cB_i,\cA_4)\\
&\geq&\prod_{i=4}^{3r-1}\brac{1-\frac{2r^7}{n^3}}\\
&\geq&1-\frac{6r^8}{n^3}.
\end{eqnarray*}
\proofend

\subsection{Proof of Lemma~\ref{propD}}

Assume that $\cB_{3r}\cap \cA_4$ occurs and that $v$ is the unique vertex of degree $3r$ in $\cF_{3r}$.
We show that \whp\ $v\in e_i$ for $i>3r$.
\begin{claim}\label{cl1}
Suppose that $\cB_{3r}\cap \cA_4$ occurs. Then
$e_i'=e_i\setminus \{v\},1\leq i\leq 3r$ is a collection of
$3r$ randomly chosen $(r-1)$-sets from $[n]\setminus\{v\}$.
\end{claim}
The claim can be argued as follows: $e_i$ is chosen uniformly from all
$r$-sets which meet $e_1,e_2,\ldots,e_{i-1}$. If we add the condition
$v\in e_i$ i.e. $\cB_i$ occurs, then $e_i$ is equally likely to be any
such $r$-set containing $v$.
\proofend

Recall that $\cD$ is the event that there is an $r$-set which meets
all edges but does not contain the kernel. Then
\begin{eqnarray*}
\Pr(\cD\mid \cB_{3r},\cA_4)&\leq&\binom{n-1}{r}
\brac{1-\frac{\binom{n-r-1}{r-1}}{\binom{n-1}{r-1}}}^{3r}\nonumber\\
&\leq&\bfrac{ne}{r}^r\bfrac{r^2}{n-2r}^{3r}\nonumber\\
&=&\bfrac{ner^5}{(n-2r)^3}^r\nonumber\\
&\leq&\left(\frac{2er^5}{n^2}\right)^r
\end{eqnarray*}

\proofend

\subsection{Proof of Lemma~\ref{propB3}}

We show that $\Pr(\cB_3)=o(1)$. We write
\begin{equation}\label{B3}
\Pr(\cB_3)=\sum_{i=1}^{r-1}f(i)g(i)
\end{equation}
where
\begin{eqnarray}
f(i)&=&\Pr(|e_1\cap e_2|=i)\nonumber\\
&=&\frac{\binom{r}{i}\binom{n-r}{r-i}}{\binom{n}{r}-\binom{n-r}{r}}\label{fi}\\
\noalign{and}
g(i)&=&\Pr(\cB_3\mid |e_1\cap e_2|=i)\nonumber\\
&=&\frac{\binom{n}{r}-\binom{n-i}{r}}{\binom{n}{r}-2\binom{n-r}{r}+\binom{n-2r+i}{r}}\label{gi}
\end{eqnarray}
Now for $0\leq s\leq 2r$ we have
\begin{eqnarray}
\frac{\binom{n-s}{r}}{\binom{n}{r}}&=&\prod_{j=0}^{r-1}\brac{1-\frac{s}{n-j}}\nonumber\\
&=&\prod_{j=0}^{r-1}\exp\left\{-\frac{s}{n}+O\bfrac{r^2}{n^2}\right\}\nonumber\\
&=&\exp\left\{-\frac{rs}{n}+O\bfrac{r^3}{n^2}\right\}.\label{rs}
\end{eqnarray}
Furthermore,
\begin{eqnarray*}
\frac{\binom{r}{i}\binom{n-r}{r-i}}{\binom{n}{r}}&\leq&\frac{r^i}{i!}\cdot\frac{\binom{n-r}{r-i}}
{\binom{n-r}{r}}\cdot\frac{\binom{n-r}{r}}{\binom{n}{r}}\\
&\leq&\frac{r^i}{i!}\cdot\frac{r^i}{(n-2r)^i}\exp\left\{-\frac{r^2}{n}+O\bfrac{r^3}{n^2}\right\}.
\end{eqnarray*}
Thus
\begin{equation}\label{fii}
f(i)\leq\frac{r^{2i}}{i!(n-2r)^i}
\cdot\frac{1+o(1)}{\exp\left\{\frac{r^2}{n}\right\}-1}.
\end{equation}
Using (\ref{rs}) in (\ref{gi}) we see that
\begin{eqnarray*}
g(i)&=&\frac{1-\exp\left\{-\frac{ir}{n}+O\bfrac{r^3}{n^2}\right\}}
{1-2\exp\left\{-\frac{r^2}{n}+O\bfrac{r^3}{n^2}\right\}+
\exp\left\{-\frac{r(2r-i)}{n}+O\bfrac{r^3}{n^2}\right\}}\\
&\leq&(1+o(1))\frac{\frac{ir}{n}+O\bfrac{r^3}{n^2}}{\brac{1-\exp\left\{-\frac{r^2}{n}\right\}}^2}.
\end{eqnarray*}

So,
\begin{eqnarray*}
\sum_{i=1}^{r-1}f(i)g(i)&\leq& (1+o(1))\frac{\exp\left\{\frac{2r^2}{n}\right\}}
{\brac{\exp\left\{\frac{r^2}{n}\right\}-1}^3}\sum_{i=1}^{r-1}
\frac{r^{2i}}{i!(n-2r)^i}\brac{\frac{ir}{n}+O\bfrac{r^3}{n^2}}\\
&=&O\brac{\frac{\exp\left\{\frac{2r^2}{n}\right\}}
{\brac{\exp\left\{\frac{r^2}{n}\right\}-1}^3}\frac{r^3}{n^2}\exp\left\{\frac{r^2}{n-2r}\right\}}\\
&=&O\brac{\frac{r^3}{n^2}\cdot\frac{1}{\brac{1-\exp\left\{-\frac{r^2}{n}\right\}}^3}}\\
&=&o(1).
\end{eqnarray*}
\proofend

\subsection{Proof of Lemma~\ref{propBm}}

Consider $m$ members of $\binom{[n]}{r}$ being chosen at random
(without replacement).

The probability that these $m$ edges fail to form an intersecting
family is at most
$$ \binom{m}{2}\frac{\binom{n-r}{r}}{\binom{n}{r}} \leq
\frac{m^2}{2}\left(1-\frac{r}{n}\right)^r \leq
\frac{m^2}{2}\exp\left\{-\frac{r^2}{n}\right\}
$$

Let us take
$$m=\exp\left\{\frac{r^2}{3n}\right\}.$$
For $r=\omega(\sqrt{n})$ we can use the fact that
$\cF_m$ has the same distribution as $m$ distinct randomly chosen
$r$-sets, conditional on the event (of probability $1-o(1)$) that
$\cF_m$ is intersecting. To see this consider sequentially choosing $m$
distinct sets at random. If we ignore the cases when the $m$ chosen sets are not intersecting
then we will produce a collection with the same distribution as $\cF_m$.

Using $r<n/2$, the probability that $\cF_m$ has a vertex of degree $m$ is at
most
\begin{eqnarray*}
\frac{1}{2}\exp\left\{-\frac{r^2}{3n}\right\}+
n\brac{\frac{\binom{n-1}{r-1}}{\binom{n}{r}}}^m
&=&O\left(\exp\left\{-\frac{r^2}{3n}\right\}\right)+r^mn^{1-m}\\
&= &O\left(\exp\left\{-\frac{r^2}{3n}\right\}\right)+n2^{-m}.
\end{eqnarray*}
\proofend

\section{Open Problem}
It is known that a maximal intersecting system, i.e, a system to which
we can not add any additional edge without making it non-intersecting,
may have various structures. Thus we finish by posing the following
problem.

\noindent{\bf Problem:} What is the structure of $\cF$ in different
ranges of $n^{1/3}\ll r<n/2$?

\end{document}